\documentclass[11pt]{article}
\usepackage[latin1]{inputenc}
\usepackage{epsfig}
\usepackage{color}
\usepackage[british,english]{babel}
\usepackage{amsthm}
\usepackage{amsmath}
\usepackage{amsfonts}
\usepackage{amssymb}
\usepackage{graphicx}
\newcommand{\mylabel}[1]{\label{#1}
            \ifx\undefined\stillediting
            \else \fbox{$#1$}\fi }
\newcommand{\BE}{\begin{equation}}

\newcommand{\EEQ}{\end{equation}}
\newcommand{\rfb}[1]{\mbox{\rm
   (\ref{#1})}\ifx\undefined\stillediting\else:\fbox{$#1$}\fi}

\newtheorem{corollary}{Corollary}[section]
\newtheorem{definition}[corollary]{Definition}

\newtheorem{lemma}[corollary]{Lemma}
\newtheorem{proposition}[corollary]{Proposition}
\newtheorem{remark}[corollary]{Remark}
\newtheorem{theorem}[corollary]{Theorem}
\newfont{\sBlackboard}{msbm10 scaled 900}

\newfont{\Blackboard}{msbm10 scaled 1200}

\newfont{\roma}{cmr10 scaled 1200}

\def\CC{\rm \hbox{C\kern-.56em\raise.4ex
         \hbox{$\scriptscriptstyle |$}\kern+0.5 em }}

\def\n{|\kern -.05cm{|}\kern -.05cm{|}}

\newcommand{\mm}    {{\hbox{\hskip 0.5pt}}}

\newcommand{\bluff} {{\hbox{\raise 15pt \hbox{\mm}}}}
\makeatletter
\def\section{\@startsection {section}{1}{\z@}{-3.5ex plus -1ex minus
    -.2ex}{2.3ex plus .2ex}{\large\bf}}
\makeatother
\def\be{\begin{equation}}
\def\ee{\end{equation}}

\date{ }
\begin{document}
\thispagestyle{empty}
\title{\bf Existence, uniqueness and homogenization of nonlinear parabolic problems with dynamical boundary conditions in perforated media }\maketitle

\author{ \center  Mar\'ia ANGUIANO\\
Departamento de An\'alisis Matem\'atico. Facultad de Matem\'aticas.\\
Universidad de Sevilla, 41012 Sevilla (Spain)\\
anguiano@us.es\\}

\vskip20pt

 \renewcommand{\abstractname} {\bf Abstract}
\begin{abstract} 
We consider a nonlinear parabolic problem with nonlinear dynamical boundary conditions of {\it pure-reactive} type in a media perforated by periodically distributed holes of size $\varepsilon$. The novelty of our work is to consider a nonlinear model where the nonlinearity also appears in the boundary. The existence and uniqueness of solution is analyzed. Moreover, passing to the limit when $\varepsilon$ goes to zero, a new nonlinear parabolic problem defined on a unified domain without holes with zero Dirichlet boundary condition and with extra-terms coming from the influence of the nonlinear dynamical boundary conditions is rigorously derived.
\end{abstract}
\bigskip\noindent

 {\small \bf AMS classification numbers:}  35B27, 35K57 \\
 
\bigskip\noindent {\small \bf Keywords:} Homogenization; perforated media; dynamical boundary conditions \newpage
\section {Introduction and setting of the problem}\label{S1}

Partial differential equations with dynamical boundary conditions have the main characteristic of involving the time derivative of the
unknown on the boundary of the domain. Its use dates back at least to 1901 in the context of heat transfer in a solid in contact with a moving fluid. From the second half of the 20th century until today, they have been studied in many disciplines (such as, diffusion phenomena in thermodynamics, phase-transition phenomena in material science, climate science, control theory and special flows in hydrodynamics).
 
Several approaches have been used for these problems in a periodically perforated domain, like Homogenization Theory. Many recent papers in the literature have treated the homogenization of elliptic problems with nonlinear boundary conditions with prescribed growth. In particular, in Cioranescu {\it et al.} \cite{Ciora_Donato_Zaki}, the authors apply the periodic unfolding method in perforated domains to a class of elliptic problems with nonlinear conditions on the boundary of the holes. The homogenization of quasilinear elliptic problems in periodically perforated domains with nonlinear Robin boundary conditions has been considered in Cabarrubias and Donato \cite{Cabarrubias}, Chourabi and Donato \cite{Chourabi1,Chourabi2} and Donato {\it et al.} \cite{Donato_Monsurro_Raimondi}.

For linear parabolic problems with linear dynamical boundary conditions of {\it pure-reactive} type in periodically perforated domains, the asymptotic behavior of the solution, when the size of the perforations tends to zero, is studied in Timofte \cite{Timofte}. But to our knowledge, there does not seem to be in the literature any study of the asymptotic behavior of the solution of nonlinear parabolic models associated to nonlinear dynamical boundary conditions of {\it pure-reactive} type in periodically perforated domains (up to the stochastic framework, see Wang and Duan \cite{Wang_Duan}).

Let us introduce the model we will be involved with in this paper. Let $\Omega$ be a bounded connected open set in $\mathbb{R}^N$ ($N\ge 2$), with smooth enough boundary $\partial \Omega$. Let $Y=[0,1]^N$ be the representative cell in $\mathbb{R}^N$ and $F$ an open subset of $Y$ with smooth enough boundary $\partial F$, such that $\bar F\subset Y$. We denote $Y^*=Y\setminus \bar F$.

For $k\in \mathbb{Z}^N$, each cell $Y_{k,\varepsilon}=\varepsilon\,k+\varepsilon\,Y$ is similar to the unit cell $Y$ rescaled to size $\varepsilon$ and $F_{k,\varepsilon}=\varepsilon\,k+\varepsilon\,F$ is similar to $F$ rescaled to size $\varepsilon$. We denote $Y^*_{k,\varepsilon}=Y_{k,\varepsilon}\setminus \bar F_{k,\varepsilon}$. We denote by $F_\varepsilon$ the set of all the holes contained in $\Omega$, i.e. $F_\varepsilon=\displaystyle\cup_{k\in K}\{F_{k,\varepsilon}:\bar{F}_{k,\varepsilon}\subset \Omega \},$
where $K:=\{k\in \mathbb{Z}^N:Y_{k,\varepsilon}\cap\Omega\ne \varnothing \}$.

Let $\Omega_\varepsilon=\Omega\backslash \bar F_\varepsilon$. By this construction, $\Omega_\varepsilon$ is a periodically perforated domain with holes of the same size as the period. 

We consider the following problem for a nonlinear reaction-diffusion equation with nonlinear dynamical boundary conditions of {\it pure-reactive} type on the surface of the holes and zero Dirichlet condition on the exterior boundary,
\begin{equation}
\left\{
\begin{array}
[c]{r@{\;}c@{\;}ll}%
\displaystyle\frac{\partial u_\varepsilon}{\partial t}-\Delta\,u_\varepsilon+\kappa u_\varepsilon+f(u_\varepsilon) &
= &
h(x,t)\quad & \text{\ in }\;\Omega_\varepsilon\times(0,T) ,\\
\displaystyle\frac{\partial
u_\varepsilon}{\partial\vec{n}}+\varepsilon\,\displaystyle\frac{\partial u_\varepsilon}{\partial
t}+\varepsilon\,g(u_\varepsilon) & = & \varepsilon\,\rho(x,t) & \text{\ on }%
\;\partial F_\varepsilon\times( 0,T),\\
u_\varepsilon(x,0) & = & u_\varepsilon^{0}(x), & \text{\ for }\;x\in\Omega_\varepsilon,\\
u_\varepsilon(x,0) & = & \psi_\varepsilon^{0}(x), & \text{\ for
}\;x\in\partial F_\varepsilon,\\
u_\varepsilon&=& 0, & \text{\ on }%
\;\partial \Omega\times( 0,T),
\end{array}
\right. \label{PDE}%
\end{equation}
where $\vec{n}$ is the outer normal to $\partial F_\varepsilon$, $T>0$, and
\begin{equation}\label{hyp 0}
\kappa>0,\quad u_\varepsilon^{0}\in L^2\left( \Omega\right),\quad
\psi_\varepsilon^{0}\in L^{2}\left( \partial F_\varepsilon\right),
\end{equation}
\begin{equation}\label{hyp 0'}
h\in L^{2}\left(0,T
;L^{2}\left( \Omega\right) \right),\quad \rho\in
L^{2}\left(0,T;H_0^{1}\left(
\Omega\right) \right),\end{equation} are given. 

We also assume that the functions $f$ and $g\in{C}\left(
\mathbb{R}\right) $ are given, and satisfy that there exist exponents $p$ and $q$ such that 
\begin{equation}\label{assumption_p}
2\leq p<+\infty,  \text{ if } N=2\quad  \text{ and } \quad 2\leq p\leq {2N \over N-2}, \text{ if } N>2,
\end{equation}
\begin{equation}\label{assumption_q}
2\leq q<+\infty,  \text{ if } N=2\quad  \text{ and } \quad 2\leq q\leq {2N \over N-2}, \text{ if } N>2,
\end{equation}
and constants $\alpha_1>0$, $\alpha_2>0,$
$\beta>0$, and $l>0$, such that
\begin{equation}
\alpha_{1}\left\vert s\right\vert ^{p}-\beta\leq f(s)s\leq\alpha
_{2}\left\vert s\right\vert ^{p}+\beta,\quad\text{for all
$s\in\mathbb{R}$,} \label{hip_1}%
\end{equation}%
\begin{equation}
\alpha_{1}\left\vert s\right\vert ^{q}-\beta\leq g(s)s\leq\alpha
_{2}\left\vert s\right\vert ^{q}+\beta,\quad\text{for all
$s\in\mathbb{R}$,} \label{hip_2}%
\end{equation}
\begin{equation}
\left( f(s)-f(r)\right) \left( s-r\right) \geq-l\left(
s-r\right) ^{2},\quad\text{for all
$s,r\in\mathbb{R}$,} \label{hip_3}%
\end{equation}and
\begin{equation}
\left( g(s)-g(r)\right) \left( s-r\right) \geq-l\left(
s-r\right) ^{2},\quad\text{for all
$s,r\in\mathbb{R}$.} \label{hip_4}%
\end{equation}
It is easy to see from (\ref{hip_1}) and (\ref{hip_2}) that there
exists a constant $C>0$
 such that%
\begin{equation}\label{hipo_consecuencia}
\left\vert f(s)\right\vert
 \leq
C\left( 1+\left\vert s\right\vert ^{p-1}\right) \text{, \ \
\ }\left\vert g(s)\right\vert \leq C\left( 1+\left\vert
s\right\vert ^{q-1}\right),\quad\text{for all
$s\in\mathbb{R}$.}
\end{equation}
Let us denote
\[
\mathcal{F}(s):=\int_{0}^{s}f(r)dr\quad\textrm{and}\quad \mathcal{G}(s):=\int_{0}^{s}g(r)dr.
\]
Then, there exist positive constants $\widetilde{\alpha}_{1}$,
$\widetilde {\alpha}_{2},$ and $\widetilde{\beta}$ such that
\begin{equation}
\widetilde{\alpha}_{1}| s|^{p}-\widetilde{\beta} \leq\mathcal{F}(s)\leq\widetilde{\alpha}_{2}|s|^{p}+\widetilde{\beta}\quad\forall s\in\mathbb{R}, \label{hip_1_adicional}
\end{equation}
and
\begin{equation}
\widetilde{\alpha}_{1}| s|^{q}-\widetilde{\beta} \leq\mathcal{G}(s)\leq\widetilde{\alpha}_{2}|s|^{q}+\widetilde{\beta}\quad\forall s\in\mathbb{R}. \label{hip_2_adicional}
\end{equation}
\begin{remark}
If $u_\varepsilon$ is regular enough, then a compatibility condition for problem (\ref{PDE}) is that $\psi_\varepsilon^{0}$ must coincide with the restriction to $\partial F_\varepsilon$ of $u_\varepsilon^0$, and therefore the fourth equation in (\ref{PDE}) is omitted. Nevertheless, this equation seems necessary for the concept of weak solution (see Definition \ref{definition_weakSolution}).
\end{remark}

In this paper, our main motivation is to study the asymptotic behavior, as $\varepsilon \to 0$, of the solution $u_\varepsilon$ of (\ref{PDE}). As we mentioned before, we only have references in the literature of this approach in the stochastic context. In that sense, a particularly interesting situation is treated in Wang and Duan \cite{Wang_Duan} with the help of the two-scale convergence. There, the authors obtain the asymptotic behavior of the solution of a stochastic partial differential equation with random dynamical boundary conditions, under the restrictive assumption $g(s)=s$, i.e. the nonlinearity does not appear in the boundary. However, we will obtain the asymptotic behavior of the solution of (\ref{PDE}) where the nonlinearity also appears in the boundary. We use the energy method of Tartar \cite{Tartar}, which has been considered by many authors (see, for instance, Cioranescu and Donato \cite{Ciora2}) and the technique introduced by Vanninathan \cite{Vanni} for the Steklov problem which transforms surface integrals into volume integrals. This technique was already used as a main tool to homogenize the non homogeneous Neumann problem for the elliptic case by Cioranescu and Donato \cite{Ciora2}.

In this sense, when $\varepsilon \to 0$, we have got a new nonlinear reaction-diffusion equation with constant coefficient defined on $\Omega\times (0,T)$, with zero Dirichlet boundary condition on the boundary, and with a constant extra-term in front of the derivative which comes from the well-balanced contribution of the dynamical part of the boundary condition on the surface of the holes.

The structure of the paper is as follows. In Section \ref{S2}, we give some notations which are used in the paper. In Section \ref{S3}, we give a weak formulation of the problem, the concept of weak solution, and establish the existence and uniqueness of solution using the monotonicity method. Some {\it a priori} estimates are rigorously stablished in Section \ref{S4}. A compactness result, which is the main key when we will pass to the limit later, is addressed in Section \ref{S5}. In Section \ref{S6}, the main goal of proving the asymptotic behavior of the solution is finally established in Theorem \ref{Main}.

\section{Some notations}\label{S2}
In this section, we give some notations which are used in the paper.

We denote by $\chi_{\Omega_\varepsilon}$ the characteristic function of the domain $\Omega_\varepsilon$.

We denote by $(\cdot,\cdot) _{\Omega_\varepsilon}$ (respectively, $(
\cdot,\cdot)_{\partial F_\varepsilon}$) the inner product in
$L^{2}(\Omega_\varepsilon)$ (respectively, in $L^{2}(\partial F_\varepsilon)$),
and by $\left\vert \cdot\right\vert _{\Omega_\varepsilon}$
(respectively, $\left\vert \cdot\right\vert
_{\partial F_\varepsilon}$) the associated norm. We also denote $(\cdot,\cdot) _{\Omega_\varepsilon}$ the inner product in $(L^2(\Omega_\varepsilon))^N$.

If $r\ne2$, we will also denote
$(\cdot,\cdot) _{\Omega_\varepsilon}$ (respectively, $(
\cdot,\cdot)_{\partial F_\varepsilon}$) the duality product between
$L^{r'}(\Omega_\varepsilon)$ and $L^{r}(\Omega_\varepsilon)$ (respectively, the duality
product between $L^{r'}(\partial F_\varepsilon)$ and
$L^{r}(\partial F_\varepsilon)$). We will denote
$|\cdot|_{r,\Omega_\varepsilon}$ (respectively $|\cdot|_{r,\partial F_\varepsilon}$)
the norm in $L^r(\Omega_\varepsilon)$ (respectively in $L^r(\partial F_\varepsilon)$).

By $\left\Vert \cdot\right\Vert _{\Omega_\varepsilon}$ we denote the norm in
$H^{1}\left(\Omega_\varepsilon\right)$, which is associated to the inner
product $((\cdot,\cdot))_{\Omega_\varepsilon}:=
(\nabla\cdot,\nabla\cdot)_{\Omega_\varepsilon}+\left(\cdot,\cdot\right)
_{\Omega_\varepsilon}.$

By $||\cdot||_{\Omega_\varepsilon,T}$ we denote the norm in $L^2(0,T;H^1(\Omega_\varepsilon))$. By $|\cdot|_{r,\Omega_\varepsilon,T}$ (respectively $|\cdot|_{r,\partial F_\varepsilon,T}$), we denote the norm in $L^r(0,T;L^r(\Omega_\varepsilon))$ (respectively $L^r(0,T;L^r(\partial F_\varepsilon))$).

We denote by $\gamma_{0}$ the trace operator $u\mapsto
u|_{\partial\Omega_\varepsilon}$. The trace operator belongs to
$\mathcal{L}(H^1(\Omega_\varepsilon), H^{1/2}(\partial\Omega_\varepsilon))$, and we will
use $\|\gamma_0\|$ to denote the norm of $\gamma_0$ in this space.

We will use $\|\cdot\|_{\partial\Omega_\varepsilon}$ to denote the
norm in $H^{1/2}(\partial\Omega_\varepsilon),$ which is given by
$\|\phi\|_{\partial\Omega_\varepsilon}=\inf\{\|v\|_{\Omega_\varepsilon}:\;
\gamma_0(v)=\phi\}$. We remember that with this norm,
$H^{1/2}(\partial\Omega_\varepsilon)$ is a Hilbert space. 

Finally, we denote by $H^r_{\partial \Omega}(\Omega_\varepsilon)$ and $H^r_{\partial \Omega}(\partial \Omega_\varepsilon)$, for $r\ge 0$, the standard Sobolev spaces which are closed subspaces of $H^r(\Omega_\varepsilon)$ and $H^r(\partial \Omega_\varepsilon)$, respectively, and the subscript $\partial \Omega$ means that, respectively, traces or functions in $\partial \Omega_\varepsilon$, vanish on this part of the boundary of $\Omega_\varepsilon$, i.e.
$$H^r_{\partial \Omega}(\Omega_\varepsilon)=\{v\in H^r(\Omega_\varepsilon):\gamma_0(v)=0  \text{ on } \partial \Omega \},$$
and
$$H^r_{\partial \Omega}(\partial \Omega_\varepsilon)=\{v\in H^r(\partial \Omega_\varepsilon):v=0  \text{ on } \partial \Omega \}.$$
Let us notice that, in fact, we can consider an element of $H^{1/2}(\partial F_\varepsilon)$ as an element of $H^{1/2}_{\partial \Omega}(\partial \Omega_\varepsilon)$.

Analogously, for $r\ge 2$, we denote
$$L^r_{\partial \Omega}(\partial \Omega_\varepsilon):=\{v\in L^r(\partial \Omega_\varepsilon):v=0  \text{ on } \partial \Omega \}.$$
Let us notice that, in fact, we can consider the given $\psi_\varepsilon^{0}$ as an element of $L^2_{\partial \Omega}(\partial \Omega_\varepsilon)$.

Let us consider the space
$$H_p:= L^{p}\left(
\Omega_\varepsilon\right) \times L_{\partial \Omega}^{p}\left( \partial\Omega_\varepsilon\right)
\text{,}\quad \forall p\ge 2,
$$with the natural inner product $ ((
v,\phi), ( w,\varphi))_{H_p}=(v,w)_{\Omega_\varepsilon}+\varepsilon
(\phi,\varphi)_{\partial F_\varepsilon},$ which in particular
induces the norm $|(\cdot,\cdot)|_{H_p}$ given by
$$|\left(
v,\phi\right)|^2_{H_p}=|v|_{\Omega_\varepsilon}^2+\varepsilon|\phi|^2_{\partial F_\varepsilon},\quad(v,\phi)\in
H_p.$$
For the sake of clarity, we shall omit to write explicitly the index $p$ if $p=2$, so we denote by $H$ the Hilbert space $$H:=L^{2}\left(
\Omega_\varepsilon\right) \times L_{\partial \Omega}^{2}\left( \partial\Omega_\varepsilon\right).$$
Let us also consider the space
$$V_1:=\left\{ \left( v,\gamma_{0}(v)\right) :v\in H_{\partial \Omega}^{1}\left(
\Omega_\varepsilon\right) \right\}.$$ We note that $V_1$ is a closed
vector subspace of $H_{\partial \Omega}^{1}\left( \Omega_\varepsilon\right) \times
H_{\partial \Omega}^{1/2}\left(\partial\Omega_\varepsilon\right),$ and therefore, with
the norm $\|(\cdot,\cdot)\|_{V_1}$ given by
\[
\left\Vert \left( v,\gamma_{0}(v)\right) \right\Vert^2
_{V_1}=\left\Vert
v\right\Vert _{\Omega_\varepsilon}^2+\left\Vert \gamma_{0}%
(v)\right\Vert^2_{\partial F_\varepsilon}, \quad\left(
v,\gamma_{0}(v)\right) \in V_1,
\]
$V_1$ is a Hilbert space.

In what follows, we shall denote by $C$ different constants which are independent of $\varepsilon$.

\section{Existence and uniqueness of solution}\label{S3}
We state in this section a result on the existence and uniqueness of solution of problem (\ref{PDE}). Instead of working directly with our equation, we will apply a general result which is a slight modification of Theorem 1.4, Chapter 2 in Lions \cite{Lions}.

In the sequel, we assume that
\begin{equation}\label{Initial_condition}
|(u_\varepsilon^0,\psi_\varepsilon^0)|_{H}\leq C.
\end{equation}

\begin{definition}\label{definition_weakSolution} A weak solution of (\ref{PDE}) is a pair of functions $(u_\varepsilon,\psi_\varepsilon)$, satisfying
\begin{equation}\label{weak0}
 u_\varepsilon\in
C([0,T];L^2(\Omega_\varepsilon)),\quad \psi_\varepsilon\in
C([0,T];L_{\partial \Omega}^2(\partial\Omega_\varepsilon)),\quad\hbox{
for all $T>0$,}
\end{equation}
\begin{equation}\label{weak1}
 u_\varepsilon\in L^2(0,T;H_{\partial \Omega}^1(\Omega_\varepsilon))\cap L^p(0,T;L^p(\Omega_\varepsilon)),
 \quad\hbox{
for all $T>0$,}
\end{equation}
\begin{equation}\label{weak2}
\psi_\varepsilon\in L^2(0,T;H_{\partial \Omega}^{1/2}(\partial\Omega_\varepsilon))\cap
L^q(0,T;L_{\partial \Omega}^q(\partial\Omega_\varepsilon)),\quad\hbox{ for all $T>0$,}
\end{equation}
\begin{equation}\label{weak3}
 \gamma_0(u_\varepsilon(t))=\psi_\varepsilon(t),\quad\hbox{ a.e. $t\in (0,T],$}
 \end{equation}
 \begin{equation}\label{weak4}
\left\{
\begin{array}{l}
 \dfrac{d}{dt}(u_\varepsilon(t),v)_{\Omega_\varepsilon}+\varepsilon\,\dfrac{d}{dt}(
\psi_\varepsilon(t),\gamma_{0}(v))_{\partial F_\varepsilon}+(\nabla u_\varepsilon(t),\nabla v)_{
\Omega_\varepsilon}+\kappa(u_\varepsilon(t),v)_{\Omega_\varepsilon}\\[2ex]
+(f(u_\varepsilon(t)),v)_{\Omega_\varepsilon}
 +\varepsilon\,(g(\psi_\varepsilon(t)),\gamma_{0}%
(v))_{\partial F_\varepsilon}
 =(h(t),v)_{\Omega_\varepsilon}
 +\varepsilon\,(\rho(t),\gamma_{0}(v))_{\partial F_\varepsilon}\\[2ex]
 \hbox{in $\mathcal{D}'(0,T)$, for all $v\in H_{\partial \Omega}^1(\Omega_\varepsilon)\cap L^p(\Omega_\varepsilon)$
 such that $\gamma_0(v)\in L_{\partial \Omega}^q(\partial \Omega_\varepsilon),$}
\end{array}
\right.
\end{equation}
\begin{equation}\label{weak5}
 u_\varepsilon(0)=u_\varepsilon^0,\quad and\quad \psi_\varepsilon(0)=\psi_\varepsilon^0.
\end{equation}
\end{definition}

\begin{remark}
In view of Theorem 3 in Chapter 5 subsection 5.9.2 in Evans \cite{Evans}, it is not difficult to prove that it $(u_\varepsilon,\psi_\varepsilon)$ satisfies
(\ref{weak1})--(\ref{weak4}), then $(u_\varepsilon,\psi_\varepsilon)$ satisfies (\ref{weak0}). The function $\psi_\varepsilon$ is the
$L_{\partial \Omega}^2(\partial\Omega_\varepsilon)$-continuous version of $\gamma_0(u_\varepsilon)$ (see
\eqref{Weak1}--\eqref{Weak3} below).
\end{remark}

We have the following result.

\begin{theorem}
\label{Existence_solution_PDE}Under the assumptions (\ref{hyp
0})--(\ref{hip_4}), there exists a unique solution
$(u_\varepsilon,\psi_\varepsilon)$
of the problem (\ref{PDE}). Moreover, this solution satisfies the
energy equality
\begin{eqnarray}\nonumber
&&\frac{1}{2}\frac{d}{dt}\left(|(u_\varepsilon(t),\psi_\varepsilon(t))|^2_{H}\right)+|\nabla
u_\varepsilon(t)|^2_{\Omega_\varepsilon} +\kappa|u_\varepsilon(t)|^2_{\Omega_\varepsilon} \\\nonumber
 &&+(f(u_\varepsilon(t)),u_\varepsilon(t))_{\Omega_\varepsilon}+\varepsilon\,(g(\psi_\varepsilon(t)),\psi_\varepsilon(t))_{\partial F_\varepsilon}\\
 &=&(h(t),u_\varepsilon(t))_{\Omega_\varepsilon}+\varepsilon\,(\rho(t),\psi_\varepsilon(t))_{\partial F_\varepsilon},\quad\mbox{a.e. $t\in(0,T).$}\label{energyequality}
\end{eqnarray}
\end{theorem}
\begin{proof}

The proof of this result is standard. For the sake of completeness, we give a sketch
of a proof.

First, we prove that $V_1$ is densely embedded in $H$. In fact, if
we consider $\left(w,\phi\right)\in H$ such that
\[
(v,w)_{\Omega_\varepsilon}+\varepsilon(\gamma_{0}(v),\phi)_{\partial F_\varepsilon}=0,\quad
\mbox{for all $v\in H_{\partial \Omega}^{1}\left( \Omega_\varepsilon\right),$}%
\]
in particular, we have
\[
(v,w)_{\Omega_\varepsilon}=0,\quad\mbox{for all $v\in H_{0}^{1}\left(
\Omega_\varepsilon\right),$}
\]
and therefore $w=0$. Consequently,
\[
(\gamma_{0}(v),\phi)_{\partial F_\varepsilon}=0,\quad\mbox{for all $v\in
H_{\partial \Omega}^{1}\left( \Omega_\varepsilon\right)$,}%
\]
and then, as $H_{\partial \Omega}^{1/2}\left( \partial\Omega_\varepsilon\right)
=\gamma_{0}\left( H_{\partial \Omega}^{1}\left( \Omega_\varepsilon\right) \right) $ is
dense in $L_{\partial \Omega}^{2}\left( \partial \Omega_\varepsilon\right),$ we have that
$\phi=0$.

\vspace{0.5cm}

Now, on the space $V_1$ we define a continuous symmetric linear
operator $A_1:V_1\rightarrow V_1^{\prime}$, given by
\begin{equation}\label{def_A1}
\langle A_1(( v,\gamma_{0}(v))) ,( w,\gamma _{0}(w))\rangle
=(\nabla v,\nabla w)_{\Omega_\varepsilon}+\kappa(v,w)_{\Omega_\varepsilon}\text{, \
}\forall v,w\in H_{\partial \Omega}^{1}\left( \Omega_\varepsilon\right)
\text{.}%
\end{equation}
We observe that $A_1$ is coercive. In fact, we have
\begin{eqnarray}\nonumber
\left\langle A_1\left( \left( v,\gamma_{0}(v)\right)
,\left( v,\gamma _{0}(v)\right) \right) \right\rangle &
\geq&\min\left\{ 1,\kappa\right\}
\left\Vert v\right\Vert _{\Omega_\varepsilon}^{2}%
\\
& = &\frac{1}{1+\|\gamma_0\|^2}\min\left\{ 1,\kappa\right\}
\left\Vert v\right\Vert
_{ \Omega_\varepsilon }^{2}\nonumber\\
&& +\frac{\|\gamma_0\|^2}{1+\left\Vert
\gamma_{0}\right\Vert ^{2}}\min\left\{
1,\kappa\right\} \left\Vert v\right\Vert _{\Omega_\varepsilon}^{2}\nonumber\\
& \geq&\frac{1}{1+\|\gamma_0\|^2}\min\left\{
1,\kappa\right\}
 \left\Vert \left( v,\gamma_{0}(v)\right) \right\Vert _{V_1}%
^{2}\text{,}\label{Coercitivity}
\end{eqnarray}for all $v\in H_{\partial \Omega}^1(\Omega_\varepsilon)$.

Let us denote
$$ V_{2}=L^{p}\left( \Omega_\varepsilon\right)
\times L_{\partial \Omega}^{2}\left( \partial\Omega_\varepsilon\right),\quad
V_{3}=L^{2}\left( \Omega_\varepsilon\right) \times L_{\partial \Omega}^{q}\left(
\partial\Omega_\varepsilon\right),$$
$$A_{2}\left( v,\phi\right) =(f(v),0),\quad
A_{3}\left( v,\phi\right) =(0,\varepsilon\,g(\phi)), \quad
\vec{h}(t)=(h(t),\varepsilon\rho(t)).
$$
From \eqref{hipo_consecuencia} one deduces that
$A_i:V_i\rightarrow V'_i$, for $i=2,3.$

Observe also that by \eqref{hyp 0'},
$$\vec{h}\in L^{2}\left(0,T;H
\right)\subset L^{2}\left(0,T;V_1^{\prime}\right).
$$

With this notation, and denoting $V=\cap_{i=1}^3V_i,$ $p_1=2,$
$p_2=p,$ $p_3=q,$ $\vec{u}_\varepsilon=(u_\varepsilon,\psi_\varepsilon)$, one has that
\eqref{weak0}--\eqref{weak5} is equivalent to
\begin{equation}\label{Weak1}
 \vec{u}_\varepsilon\in C([0,T];H),\quad \vec{u}_\varepsilon\in\bigcap_{i=1}^3 L^{p_i}(0,T;V_i),
 \quad\mbox{for all $T>0,$}
\end{equation}
\begin{equation}\label{Weak2}
(\vec{u}_\varepsilon)'(t)+\sum_{i=1}^3A_i(\vec{u}_\varepsilon(t))=\vec{h}(t)\quad\mbox{in
$\mathcal{D}'(0,T;V'),$}
\end{equation}
\begin{equation}\label{Weak3}
\vec{u}_\varepsilon(0)=(u_\varepsilon^0,\psi_\varepsilon^0).
\end{equation}

Applying a slight modification of
\cite[Ch.2,Th.1.4]{Lions}, it is not difficult to see that problem
\eqref{Weak1}--\eqref{Weak3} has a unique solution. Moreover,
 $\vec{u}_\varepsilon$ satisfies the energy equality
$$\frac{1}{2}\frac{d}{dt}|\vec{u}_\varepsilon(t)|^2_H+\sum_{i=1}^3\langle A_i(\vec{u}_\varepsilon(t)),\vec{u}_\varepsilon(t)\rangle_i=
(\vec{h}(t),\vec{u}_\varepsilon(t))_H\quad\mbox{a.e. $t\in(0,T),$}$$ where
 $\left\langle \cdot ,\cdot \right\rangle_i $ denotes the duality product
between $V_i^{\prime}$ and $V_i$.

This last equality turns out to be just \eqref{energyequality}.
\end{proof}
\begin{remark}
The assumption $\kappa>0$ is not necessary for the existence and uniqueness of weak solution to (\ref{PDE}).
\end{remark}

\section{A priori estimates}\label{S4}
Let us begin with a variant of the Trace Theorem in $\Omega_\varepsilon$.
\begin{lemma}\label{lemma_Evans}
There exists a positive constant $C$ independent of $\varepsilon$, such that
\begin{equation}\label{Evans}
\varepsilon|\gamma_0(v)|^p_{p,\partial F_\varepsilon}\leq C\left(|v|^p_{p,\Omega_\varepsilon}+\varepsilon^p|\nabla v|^p_{p,\Omega_\varepsilon} \right), \quad 1\leq p<\infty,
\end{equation}
for any $v\in W^{1,p}(\Omega_\varepsilon)$, $v=0$ on $\partial \Omega$.
\end{lemma}
\begin{proof}
For any function $v(y)\in W^{1,p}(Y^*)^N$, using the Trace Theorem (see Chapter 5, Section 5.5, Theorem 1 in Evans \cite{Evans}, for more details), we have for every $k\in \mathbb{Z}^N$
\begin{equation}\label{Trace}
\int_{\partial (F+k)}|\gamma_0(v)|^pd\sigma(y)\leq C\left(\int_{Y^*+k} |v|^pdy+\int_{Y^*+k}|\nabla v|^pdy\right), \quad 1\leq p<\infty,
\end{equation}
where the constant $C$ depends only on $p$ and $Y$.

By the change of variable
\begin{equation}\label{dilatacion}
y={x \over \varepsilon},\quad d\sigma(y)=\varepsilon^{-(N-1)}d\sigma(x),\quad \partial_y=\varepsilon\, \partial_x,
\end{equation}
we rescale (\ref{Trace}) from $Y^*+k$ to $Y^*_{k,\varepsilon}$ and from $F+k$ to $F_{k,\varepsilon}$. This yields that, for any function $v(x)\in W^{1,p}(Y^*_{k,\varepsilon})^N$, one has
\begin{equation*}
\varepsilon\int_{\partial F_{k,\varepsilon}}|\gamma_0(v)|^pd\sigma(x)\leq C\left(\int_{Y^*_{k,\varepsilon}}|v|^pdx+\varepsilon^p\int_{Y^*_{k,\varepsilon}}|\nabla v|^pdx \right),
\end{equation*}
with the same constant $C$ as in (\ref{Trace}). Summing the inequalities, for every $k\in K$, gives the desired result (\ref{Evans}).

\end{proof}

Let us obtain some {\it a priori} estimates for $u_\varepsilon$.
\begin{lemma}\label{estimates1}
Under the assumptions (\ref{hyp 0})--(\ref{hip_4}) and (\ref{Initial_condition}), there exists a constant $C$ independent of $\varepsilon$, such that the solution $u_\varepsilon$ of the problem (\ref{PDE}) satisfies
\begin{equation}\label{acotacion2}
\left\Vert
u_\varepsilon\right\Vert _{\Omega_\varepsilon,T}\leq C.
\end{equation}
\end{lemma}
\begin{proof}
By \eqref{energyequality} and
taking into account \eqref{hip_1}, \eqref{hip_2} and (\ref{Coercitivity}), we have
\begin{eqnarray}\nonumber
&&\frac{d}{dt}\left(|(u_\varepsilon(t),\psi_\varepsilon(t))|^2_{H}\right)+\frac{2\min\left\{
1,\kappa\right\}}{1+\|\gamma_0\|^2}\left(\left\Vert
u_\varepsilon(t)\right\Vert _{\Omega_\varepsilon}^2+\left\Vert \psi_\varepsilon(t)\right\Vert^2_{\partial F_\varepsilon}\right) \\\nonumber&&  +2\alpha_1(|u_\varepsilon(t)|_{p,\Omega_\varepsilon}^p+\varepsilon\,|\psi_\varepsilon(t)|_{q,\partial F_\varepsilon}^q)
 \\ \label{precont1}
 &\leq&
 2\beta(|\Omega_\varepsilon|+\varepsilon\,|\partial F_\varepsilon|)
 +|h(t)|^2_{\Omega_\varepsilon}+\varepsilon\,|\rho(t)|^2_{\partial F_\varepsilon}
 +|u_\varepsilon(t)|^2_{\Omega_\varepsilon}+\varepsilon\,|\psi_\varepsilon(t)|^2_{\partial F_\varepsilon}, 
\end{eqnarray}
where $|\Omega_\varepsilon|$ and $|\partial F_\varepsilon|$ denote the measure of $\Omega_\varepsilon$ and $\partial F_\varepsilon$, respectively.

Integrating (\ref{precont1}) between $0$ and $t$ and using (\ref{hyp 0'}), we obtain
\begin{eqnarray}\nonumber
&&|(u_\varepsilon(t),\psi_\varepsilon(t))|^2_{H}+\frac{2\min\left\{
1,\kappa\right\}}{1+\|\gamma_0\|^2}\int_0^t\left(\left\Vert
u_\varepsilon(s)\right\Vert _{\Omega_\varepsilon}^2+\left\Vert \psi_\varepsilon(s)\right\Vert^2_{\partial F_\varepsilon}\right)ds\\\nonumber&&  +2\alpha_1\int_0^t(|u_\varepsilon(s)|_{p,\Omega_\varepsilon}^p+\varepsilon\,|\psi_\varepsilon(s)|_{q,\partial F_\varepsilon}^q)ds
 \\ \label{beforeintegrating2}
 &\leq&
 2\beta t(|\Omega_\varepsilon|+\varepsilon\,|\partial F_\varepsilon|)
+ |(u_\varepsilon^0,\psi_\varepsilon^0)|^2_{H}+\int_0^T\left(|h(s)|^2_{\Omega_\varepsilon}+\varepsilon\,|\rho(s)|^2_{\partial F_\varepsilon}
\right)ds\\&&
 +\int_0^t|(u_\varepsilon(s),\psi_\varepsilon(s))|^2_{H}ds.\nonumber
\end{eqnarray}

By Lemma \ref{lemma_Evans} with $p=2$, we can deduce
\begin{equation*}
\varepsilon\,|\rho(t)|^2_{\partial F_\varepsilon}\leq C\left(|\rho(t)|^2_{\Omega_\varepsilon}+\varepsilon^2 |\nabla \rho(t)|^2_{\Omega_\varepsilon}\right)\leq C ||\rho(t)||^2_{\Omega_\varepsilon},
\end{equation*}
which together with (\ref{hyp 0'}) gives
\begin{equation}\label{trace1}
\int_0^T\left(|h(s)|^2_{\Omega_\varepsilon}+\varepsilon\,|\rho(s)|^2_{\partial F_\varepsilon}
\right)ds\leq C.
\end{equation}

On the other hand, the number of holes is given by $$N(\varepsilon)={|\Omega| \over (2\varepsilon)^N}\left(1+o(1)\right),$$
then using the change of variable (\ref{dilatacion}), we can deduce
$$|\partial F_\varepsilon|=N(\varepsilon) |\partial F_{k,\varepsilon}|=N(\varepsilon)\varepsilon^{N-1}|\partial F|\leq {C \over \varepsilon}.$$

And since $|\Omega_\varepsilon|\leq |\Omega|$, we have that 
\begin{equation}\label{trace2}
|\Omega_\varepsilon|+\varepsilon\,|\partial F_\varepsilon|\leq C.
\end{equation}

Taking into account (\ref{Initial_condition}), (\ref{trace1})-(\ref{trace2}) in (\ref{beforeintegrating2}) and applying Gronwall Lemma, in particular we obtain that there exists a positive constant $C$ such that
\begin{eqnarray}\label{despues_gronwall}
|(u_\varepsilon(t),\psi_\varepsilon(t))|^2_{H} \leq C,
\end{eqnarray}
for all $t\in (0,T)$.

Now, taking into account (\ref{Initial_condition}), (\ref{trace1})-(\ref{despues_gronwall}) in (\ref{beforeintegrating2}), we get (\ref{acotacion2}).
\end{proof}

Now, if we want to take the inner product in (\ref{PDE}) with $u'_\varepsilon$, we need that $u'_\varepsilon \in L^2(0,T;H_{\partial \Omega}^1(\Omega_\varepsilon))\cap L^p(0,T;L^p(\Omega_\varepsilon))$ with $\gamma_0(u'_\varepsilon)\in L^q(0,T;L^q_{\partial \Omega}(\Omega_\varepsilon))$. However, we do not have it for our weak solution. Therefore, we use the Galerkin method and the following lemma in order to prove, rigorously, new {\it a priori} estimates for $u_\varepsilon$.

\begin{lemma}[Lemma 11.2 in Robinson \cite{Robinson}]
\label{lema_pepe}Let $X,Y$ be Banach spaces such that $X$ is reflexive, and the inclusion $X\subset Y$ is compact.
Assume that $\{u_{m}\}$ is uniformly bounded in $L^{\infty}(0,T;X)$,
$${\rm ess}\sup_{t\in[0,T]}\|u_m(t)\|_{X}\leq C,$$
and that $u_{m} \rightharpoonup u$ weakly in $L^{2}(0,T;X)$, then
$${\rm ess}\sup_{t\in[0,T]}\|u(t)\|_{X}\leq C.$$
Furthermore, if $u\in C([0,T];Y)$, then $u(t)\in X$ for all $t\in [0,T]$ and
\[
\sup_{t\in[0,T]}\|u(t)\|_{X}\leq C.
\]
\end{lemma}

Let us observe that the space $H^{1}_{\partial \Omega}(\Omega_\varepsilon)\times H_{\partial \Omega}^{1/2}(\partial\Omega_\varepsilon)$ is compactly imbedded in $H$, and
therefore, for the symmetric and coercive linear continuous operator $A_{1}:V_{1}\rightarrow V_{1}^{\prime}$, where $A_1$ is given by (\ref{def_A1}), there exists a non-decreasing sequence $0<\lambda_{1}\leq\lambda_{2}\leq\ldots$ of eigenvalues associated to the operator
$A_{1}$ with $\lim_{j\rightarrow\infty }\lambda_{j}=\infty,$ and there exists a Hilbert basis of $H$, $\{(w_{j},\gamma_{0}(w_{j})) :j\geq1\}$$\subset D(A_1)$, with $span\{(w_{j},\gamma_{0}(w_{j})):j\geq1\} $ densely embedded in $V_{1}$, such that
\[
A_{1}((w_{j},\gamma_{0}(w_{j})))=\lambda_{j}(w_{j},\gamma_{0}(w_{j}))\quad\forall j\geq1.
\]

\begin{remark}\label{remark_densidad}
It can be proved that if $p=q\ge 2$, then $span\{(w_{j},\gamma_{0}(w_{j})):j\geq1\} $ is densely embedded in $V_1\cap H_p$.
\end{remark}

Taking into account the above facts, we denote by $$(u_{\varepsilon,m}(t),\gamma_{0}(u_{\varepsilon,m}(t)))=(u_{\varepsilon,m}(t;0,u_\varepsilon^0,\psi_\varepsilon^0),\gamma_{0}(u_{\varepsilon,m}(t;0,u_\varepsilon^0,\psi_\varepsilon^0)))$$ the Galerkin approximation of the solution
$(u_\varepsilon(t;0,u_\varepsilon^0,\psi_\varepsilon^0),\gamma_{0}(u_\varepsilon(t;0,u_\varepsilon^0,\psi_\varepsilon^0)))$ to (\ref{PDE}) for each integer $m\geq1$, which is given by
\begin{equation}
(u_{\varepsilon,m}(t),\gamma_{0}(u_{\varepsilon,m}(t)))=\sum_{j=1}^{m}\delta_{\varepsilon mj}(t)(w_{j},\gamma_{0}(w_{j})),\label{Galerkin1}
\end{equation}
and is the solution of
\begin{eqnarray}
\nonumber &&\dfrac{d}{dt}((u_{\varepsilon,m}(t),\gamma_{0}(u_{\varepsilon,m}(t))),(w_{j},\gamma_{0}(w_{j})))_{H}+\left\langle A_{1}((u_{\varepsilon,m}(t),\gamma_{0}(u_{\varepsilon,m}(t)))),(w_{j},\gamma_{0}(w_{j}))\right\rangle\\
\nonumber &&+(f(u_{\varepsilon,m}(t)),w_{j})_{\Omega_\varepsilon}+\varepsilon(  g(\gamma_{0}(u_{\varepsilon,m}(t))),\gamma_{0}(w_{j}))_{\partial F_\varepsilon}\\
&=&(  h(t),w_{j})_{\Omega_\varepsilon}+\varepsilon(  \rho(t),\gamma_{0}(w_{j}))_{\partial F_\varepsilon},\quad j=1,\ldots ,m,\label{7}
\end{eqnarray}
with initial data
\begin{equation}
\label{7'}
(u_{\varepsilon,m}(0),\gamma_{0}(u_{\varepsilon,m}(0)))=(u_{\varepsilon,m}^{0},\gamma_{0}(u_{\varepsilon,m}^{0})),
\end{equation}
where
\[
\delta_{\varepsilon mj}(t)=(u_{\varepsilon,m}(t),w_{j})_{\Omega_\varepsilon}+( \gamma_{0}(u_{\varepsilon,m}(t)),\gamma_{0}(w_{j}))_{\partial F_\varepsilon},
\]
and $(u_{\varepsilon,m}^{0},\gamma_0(u_{\varepsilon,m}^{0}))\in span\{(w_j,\gamma_0(w_j)): j=1,\ldots ,m\}$ converge (when
$m\to\infty$) to $(u_\varepsilon^0,\psi_\varepsilon^0)$ in a suitable sense which will be specified below.\\

\begin{lemma}\label{estimates2}
Suppose that in addition to the assumptions (\ref{hyp 0})--(\ref{hip_4}), we have $p=q\ge 2$. Then, for any initial condition $(u_\varepsilon^0,\psi_\varepsilon^0)\in V_1\cap H_p$, there exists a constant $C$ independent of $\varepsilon$, such that the solution $u_\varepsilon$ of the problem (\ref{PDE}) satisfies
\begin{equation}\label{acotacion5}
 \sup_{t\in [0,T]}\left\Vert u_\varepsilon(t)\right\Vert _{\Omega_\varepsilon}\leq C.
\end{equation}
\end{lemma}
\begin{proof}
Let $(u_\varepsilon^0,\psi_\varepsilon^0)\in V_{1}\cap H_p$. For all $m\geq1$, by Remark
\ref{remark_densidad}, there exists $(u_{\varepsilon,m}^{0},\gamma_{0}(u_{\varepsilon,m}^{0}))\in span\{(w_{j},\gamma_{0}(w_{j})):1\leq
j\leq m\} $, such that $\{(u_{\varepsilon,m}^{0},\gamma_{0}(u_{\varepsilon,m}^{0}))\}$ converges to $(u_\varepsilon^0,\psi_\varepsilon^0)$ in $V_{1}$
and in $H_p$. Then, in particular we know that there exists a constant $C$ such that
\begin{equation}\label{Initial_condition2}
||(u_{\varepsilon,m}^{0},\gamma_{0}(u_{\varepsilon,m}^{0}))||_{V_1}\leq C,\quad | (u_{\varepsilon,m}^{0},\gamma_{0}(u_{\varepsilon,m}^{0})) |_{H_p}\leq C.
\end{equation}
For each integer $m\geq1$, we consider the sequence $\{(u_{\varepsilon,m}(t),\gamma_{0}(u_{\varepsilon,m}(t)))\}$ defined by
(\ref{Galerkin1})-(\ref{7'}) with these initial data.

Multiplying by the derivative $\delta'_{\varepsilon mj}$ in (\ref{7}), and summing from $j=1$ to $m$, we obtain
\begin{eqnarray*}
&& |(u_{\varepsilon,m}^{\prime}(t),\gamma_{0}(u_{\varepsilon,m}^{\prime}(t)))|^2_{H}+\frac{1}{2}\frac{d}{dt}(\left\langle A_{1}((u_{\varepsilon,m}(t),\gamma_{0}(u_{\varepsilon,m}(t)))),(u_{\varepsilon,m}(t),\gamma_{0}(u_{\varepsilon,m}(t)))\right\rangle)\\
&&  +(f(u_{\varepsilon,m}(t)),u_{\varepsilon,m}^{\prime}(t))_{\Omega_\varepsilon}+\varepsilon(g(\gamma_{0}(u_{\varepsilon,m}(t))),\gamma_{0}(u_{\varepsilon,m}^{\prime}(t)))_{\partial F_\varepsilon}\\
&  =&(h(t),u_{\varepsilon,m}^{\prime}(t))_{\Omega_\varepsilon}+\varepsilon(\rho (t),\gamma_{0}(u_{\varepsilon,m}^{\prime}(t)))_{\partial F_\varepsilon}.
\end{eqnarray*}
We observe that
\[
(f(u_{\varepsilon,m}(t)),u_{\varepsilon,m}^{\prime}(t))_{\Omega_\varepsilon}=\frac{d}{dt}\int_{\Omega_\varepsilon}\mathcal{F} (u_{\varepsilon,m}(t))dx,
\]
and
$$
(g(\gamma_{0}(u_{\varepsilon,m}(t))),\gamma_{0}(u_{\varepsilon,m}^{\prime}(t)))_{\partial F_\varepsilon}
=\frac{d}{dt}\int_{\partial F_\varepsilon}\mathcal{G}(\gamma_{0}(u_{\varepsilon,m}(t))) d\sigma(x).
$$
Then, we deduce
\begin{eqnarray*}\nonumber
&&|(u_{\varepsilon,m}^{\prime}(t),\gamma_{0}(u_{\varepsilon,m}^{\prime}(t)))|^2_{H}+\frac{1}{2}\frac{d}{dt}(\left\langle A_{1}((u_{\varepsilon,m}(t),\gamma_{0}(u_{\varepsilon,m}(t)))),(u_{\varepsilon,m}(t),\gamma_{0}(u_{\varepsilon,m}(t)))\right\rangle)\\
&  \leq& \frac{1}{2}|h(t)|_{\Omega_\varepsilon}^{2}+\frac{1}{2}|u_{\varepsilon,m}^{\prime}(t)|_{\Omega_\varepsilon}^{2}+\varepsilon\frac{1}{2}|\rho(t)|_{\partial F_\varepsilon}^{2} +\varepsilon\frac{1}{2}|\gamma_{0}(u_{\varepsilon,m}^{\prime}(t))|_{\partial F_\varepsilon}^{2}\nonumber\\
&&  -\frac{d}{dt}\int_{\Omega_\varepsilon}\mathcal{F}(u_{\varepsilon,m}(t))dx-\varepsilon\frac{d}{dt}
\int_{\partial F_\varepsilon}\mathcal{G}(\gamma_{0}(u_{\varepsilon,m}(t)))d\sigma(x).\label{acotacion_necesaria}
\end{eqnarray*}
Integrating now between $0$ and $t$, taking into account the definition of $A_1$ and (\ref{Coercitivity}), and using (\ref{hip_1_adicional}) and (\ref{hip_2_adicional}), we obtain that
\begin{eqnarray}\label{last_estimate}
&& \int_{0}^{t}|(u_{\varepsilon,m}^{\prime}(s),\gamma_{0}(u_{\varepsilon,m}^{\prime}(s)))|^2_{H}ds+ \frac {\min\{1,\kappa\}}{1+\|\gamma_{0}\|^{2}}\|(u_{\varepsilon,m}(t),\gamma_{0}(u_{\varepsilon,m}(t)))\|_{V_{1}}^{2}\nonumber\\
&&
+2\widetilde{\alpha}_{1}|(u_{\varepsilon,m}(t),\gamma_{0}(u_{\varepsilon,m}(t)))|_{H_p}^{p} \nonumber\\
&\leq&\max\{1,\kappa\}\|(u_{\varepsilon,m}^{0},\gamma_{0}(u_{\varepsilon,m}^{0}))\|_{V_{1}}^{2}+\int_{0}^{T}(|h(s)|_{\Omega_\varepsilon}^{2}+\varepsilon |\rho(s)|_{\partial F_\varepsilon}^{2})ds\nonumber\\
&&
+2\widetilde {\alpha}_{2}| (u_{\varepsilon,m}^{0}, \gamma_{0}(u_{\varepsilon,m}^{0}))|_{H_p}^{p} +4\tilde\beta (|\Omega_\varepsilon|+\varepsilon|\partial F_\varepsilon|),
\end{eqnarray}
for all $t\in (0,T)$.

Taking into account (\ref{Initial_condition2}) and using (\ref{trace1})-(\ref{trace2}) in (\ref{last_estimate}), we have proved that the sequence $\{(u_{\varepsilon,m},\gamma_0(u_{\varepsilon,m}))\}$ is bounded in $C([0,T];V_1\cap H_p),$  and $\{(u_{\varepsilon,m}',\gamma_0(u_{\varepsilon,m}'))\}$ is bounded in $L^2(0, T;H),$
for all $T>0$.

If we work with the truncated Galerkin equations (\ref{Galerkin1})-(\ref{7'}) instead of the full PDE, we note that the calculations of the proof of Lemma \ref{estimates1} can be following identically to show that $\{(u_{\varepsilon,m},\gamma_0(u_{\varepsilon,m}))\}$ is bounded in $L^2(0,T;V_1),$
for all $T>0$.

Moreover, taking into account the uniqueness of solution to (\ref{PDE}) and using Aubin-Lions compactness lemma (e.g., cf. Lions \cite{Lions}), it is not difficult to conclude that the sequence $\{( u_{\varepsilon,m},\gamma_{0} (u_{\varepsilon, m}))\}$ converges weakly in $
L^{2}(0,T;V_1)$ to the solution
$(u_\varepsilon,\gamma_{0}(u_\varepsilon))$ to (\ref{PDE}). Since the inclusion $H^1(\Omega_\varepsilon)\subset L^2(\Omega_\varepsilon)$ is compact and $u_\varepsilon\in C([0,T];L^2(\Omega_\varepsilon))$, it follows using Lemma \ref{lema_pepe} that the estimate (\ref{acotacion5}) is proved.

\end{proof}

\begin{lemma}\label{estimates3}
Under the assumptions in Lemma \ref{estimates2}, assume that $f$, $g\in\mathcal{C}^{1}( \mathbb{R})$, $h\in W^{1,2}(0,T;L^{2}(\Omega))$ and $\rho\in W^{1,2}(0,T;H_0^{1}(\Omega))$, then there exists a constant $C$ independent of $\varepsilon$, such that the solution $u_\varepsilon$ of the problem (\ref{PDE}) satisfies
\begin{equation}\label{acotacion8}
\left\Vert
u'_\varepsilon\right\Vert _{\Omega_\varepsilon,T}\leq C.
\end{equation}
\end{lemma}
\begin{proof}
We first note that under the conditions imposed we have that
\begin{equation}
f^{\prime}(s)\geq-l,\qquad g^{\prime}(s)\geq-l\quad\forall s\in\mathbb{R}. \label{47}
\end{equation}
As we are assuming that $f$, $g\in\mathcal{C}^{1}(\mathbb{R})$, $h\in W_{loc}^{1,2}(\mathbb{R};L^{2}(\Omega))$ and $\rho\in W_{loc}^{1,2}(\mathbb{R};H_0^{1}(\Omega))$, we can differentiate with respect to time in (\ref{7}), and then, multiplying by the derivative $\delta'_{\varepsilon mj}$ and summing from $j=1$ to $m$, we obtain
\begin{eqnarray*}
&&  \frac{1}{2}\frac{d}{dt}|(u_{\varepsilon,m}^{\prime}(t),\gamma_{0}(u_{\varepsilon,m}^{\prime}(t)))|_{H}^{2}+\left\langle A_{1}((u_{\varepsilon,m}^{\prime}(t),\gamma_{0}(u_{\varepsilon,m}^{\prime}(t)))),(u_{\varepsilon,m}^{\prime}(t),\gamma_{0}(u_{\varepsilon,m}^{\prime}(t)))\right\rangle \\
&=&-(f^{\prime}(u_{\varepsilon,m}(t))u_{\varepsilon,m}^{\prime}(t),u_{\varepsilon,m}^{\prime}(t))_{\Omega_\varepsilon}-\varepsilon(g^{\prime}(\gamma_{0}(u_{\varepsilon,m}(t))) \gamma_{0}(u_{\varepsilon,m}^{\prime}(t)),\gamma_{0}(u_{\varepsilon,m}^{\prime}(t)))_{\partial F_\varepsilon}\\
& &+(h^{\prime}(t),u_{\varepsilon,m}^{\prime}(t))_{\Omega_\varepsilon}+\varepsilon( \rho^{\prime}(t),\gamma_{0}(u_{\varepsilon,m}^{\prime}(t)))_{\partial F_\varepsilon}.
\end{eqnarray*}
Then, using (\ref{Coercitivity}) and (\ref{47}), we have
\begin{eqnarray*}
&&\frac{d}{dt}|(u_{\varepsilon,m}^{\prime}(t),\gamma_{0}(u_{\varepsilon,m}^{\prime}(t)))|_{H}^{2}+\frac{2\min\left\{
1,\kappa\right\}}{1+\|\gamma_0\|^2}\left\Vert
(u'_{\varepsilon,m}(t),\gamma_{0}(u'_{\varepsilon,m}(t)))\right\Vert^2_{V_1}\\ \!\!&\!\leq\!&\! (2l+1)|(u_{\varepsilon,m}^{\prime}(t),\gamma_{0} (u_{\varepsilon,m}^{\prime}(r)))|_{H}^{2} +|h^{\prime}(t)|_{\Omega_\varepsilon}^{2}+\varepsilon|\rho^{\prime}(t)|_{\partial F_\varepsilon}^{2}.
\end{eqnarray*}
Integrating between $r$ and $t$
\begin{eqnarray*}
&&|(u_{\varepsilon,m}^{\prime}(t),\gamma_{0}(u_{\varepsilon,m}^{\prime}(t)))|_{H}^{2}+\frac{2\min\left\{
1,\kappa\right\}}{1+\|\gamma_0\|^2}\int_r^t \left\Vert
(u'_{\varepsilon,m}(s),\gamma_{0}(u'_{\varepsilon,m}(s)))\right\Vert^2_{V_1}ds \nonumber\\
 &\!\!  \leq\!\!&\!\!| (  u_{\varepsilon,m}^{\prime}(r),\gamma_{0}(u_{\varepsilon,m}^{\prime}(r)))  |_{H}^{2} \!+\!(2l+1)\int_{r}^{t}|(u_{\varepsilon,m}^{\prime}(s),\gamma_{0}(u_{\varepsilon,m}^{\prime}(s)))|_{H}^{2}ds \\ && +\int_{r}^{t}(|h^{\prime}(s)|_{\Omega_\varepsilon}^{2}+\varepsilon| \rho^{\prime }(s)|_{\partial F_\varepsilon}^{2})ds,\nonumber
\end{eqnarray*}
for all $0\leq r\leq t$. Now, integrating with respect to $r$ between $0$ and $t$,
\begin{eqnarray}\label{last_estimate2}
&&t|(u_{\varepsilon,m}^{\prime}(t),\gamma_{0}(u_{\varepsilon,m}^{\prime}(t)))|_{H}^{2}+\frac{2\min\left\{
1,\kappa\right\}}{1+\|\gamma_0\|^2}\int_0^t\left\Vert
(u'_{\varepsilon,m}(s), \gamma_{0}(u'_{\varepsilon,m}(s)))\right\Vert^2_{V_1}ds \\
 &\!\!  \leq\!\!& 2(l+1)\int_{0}^{t}|(u_{\varepsilon,m}^{\prime}(s),\gamma_{0}(u_{\varepsilon,m}^{\prime}(s)))|_{H}^{2}ds  +\int_{0}^{T}(|h^{\prime}(s)|_{\Omega_\varepsilon}^{2}+\varepsilon| \rho^{\prime }(s)|_{\partial F_\varepsilon}^{2})ds,\nonumber
\end{eqnarray}
for all $t\in(0,T)$.

By Lemma \ref{lemma_Evans} with $p=2$, we can deduce
\begin{equation*}
\varepsilon\,|\rho'(t)|^2_{\partial F_\varepsilon}\leq C\left(|\rho'(t)|^2_{\Omega_\varepsilon}+\varepsilon^2 |\nabla \rho'(t)|^2_{\Omega_\varepsilon}\right)\leq C ||\rho'(t)||^2_{\Omega_\varepsilon},
\end{equation*}
which, taking into account that $h\in W^{1,2}(0,T;L^{2}(\Omega))$ and $\rho\in W^{1,2}(0,T;H_0^{1}(\Omega))$, gives
\begin{equation}\label{trace1_1}
\int_0^T\left(|h'(s)|^2_{\Omega_\varepsilon}+\varepsilon\,|\rho'(s)|^2_{\partial F_\varepsilon}
\right)ds\leq C.
\end{equation}
In particular, taking into account (\ref{trace1})-(\ref{trace2}) in (\ref{last_estimate}), we have 
$$\int_0^t\left\vert
(u'_{\varepsilon,m}(s), \gamma_{0}(u'_{\varepsilon,m}(s)))\right\vert^2_{H}ds 
\leq C\left(1+\|(u_{\varepsilon,m}^{0},\gamma_{0}(u_{\varepsilon,m}^{0}))\|_{V_{1}}^{2}+| (u_{\varepsilon,m}^{0}, \gamma_{0}(u_{\varepsilon,m}^{0}))|_{H_p}^{p} \right),$$
which, jointly with (\ref{last_estimate2})-(\ref{trace1_1}), yields that
\begin{eqnarray}\label{lower1}
\int_0^t\left\Vert
(u'_{\varepsilon,m}(s), \gamma_{0}(u'_{\varepsilon,m}(s)))\right\Vert^2_{V_1}ds 
\leq C\left(1+\|(u_{\varepsilon,m}^{0},\gamma_{0}(u_{\varepsilon,m}^{0}))\|_{V_{1}}^{2}+| (u_{\varepsilon,m}^{0}, \gamma_{0}(u_{\varepsilon,m}^{0}))|_{H_p}^{p} \right),
\end{eqnarray}
and using (\ref{Initial_condition2}) we have proved that the sequence $\{(u_{\varepsilon,m}',\gamma_0(u_{\varepsilon,m}'))\}$ is bounded in $L^2(0, T;V_1),$ for all $T>0$. Then, the sequence $\{( u'_{\varepsilon,m},\gamma_{0} (u'_{\varepsilon, m}))\}$ converges weakly in $
L^{2}(0,T;V_1)$ to $(u'_\varepsilon,\gamma_{0}(u'_\varepsilon))$, for all $T>0$, and using the lower-semicontinuity of the norm and (\ref{lower1}), in particular we get
\begin{eqnarray*}
||u'_\varepsilon||_{\Omega_\varepsilon,T}&\leq& \liminf_{m\to \infty}||u'_{\varepsilon,m}||_{\Omega_\varepsilon,T}\\
&\leq& C\liminf_{m\to \infty}\left(1+\|(u_{\varepsilon,m}^{0},\gamma_{0}(u_{\varepsilon,m}^{0}))\|_{V_{1}}^{2}+| (u_{\varepsilon,m}^{0}, \gamma_{0}(u_{\varepsilon,m}^{0}))|_{H_p}^{p} \right) \\
&=&C\left(1+\|(u_\varepsilon^0,\psi_\varepsilon^0)\|_{V_{1}}^{2}+| (u_\varepsilon^0,\psi_\varepsilon^0)|_{H_p}^{p} \right),
\end{eqnarray*}
which, jointly with $(u_\varepsilon^0,\psi_\varepsilon^0)\in V_1\cap H_p$, implies (\ref{acotacion8}).

\end{proof}

\subsection{The extension of $u_\varepsilon$ to the whole $\Omega \times (0,T)$}
Since the solution $u_\varepsilon$ of the problem (\ref{PDE}) is defined only in $\Omega_\varepsilon\times (0,T)$, we need to extend it to the whole $\Omega\times (0,T)$. We denote by $\tilde v$ the extension to the whole $\Omega\times (0,T)$ for any function $v$ defined on $\Omega_\varepsilon\times (0,T)$. For finding a suitable extension $\tilde u_\varepsilon$ into all $\Omega\times (0,T)$, we shall use the following well-known extension Lemma. 
\begin{lemma}[Lemma 1 in Cioranescu and Saint Jean Paulin \cite{Cioranescu}]\label{extensionCiora}
Every function $\varphi_\varepsilon\in H^1(\Omega_\varepsilon)$, with $\varphi_\varepsilon=0$ on $\partial \Omega$, can be extended to a function $\tilde \varphi_\varepsilon\in H_0^1(\Omega)$, such that $$|\nabla \tilde \varphi_\varepsilon|_{\Omega}\leq C|\nabla \varphi_\varepsilon|_{\Omega_\varepsilon},$$
where the constant $C$ does not depend on $\varepsilon$.
\end{lemma}
Let us obtain some {\it a priori} estimates for the extension of $u_\varepsilon$ to the whole $\Omega \times (0,T)$.
\begin{corollary}\label{estimates_extension}
Assume the assumptions in Lemma \ref{estimates3}. Then, there exists an extension $\tilde u_\varepsilon$ of the solution $u_\varepsilon$ of the problem (\ref{PDE}) into $\Omega\times (0,T)$, such that
\begin{equation}\label{acotacion1_extension}
\left\Vert
\tilde u_\varepsilon(t)\right\Vert _{\Omega,T}\leq C, \quad |\tilde u_\varepsilon|_{p,\Omega,T}\leq C,
\end{equation}
\begin{equation}\label{acotacion2_extension}
\sup_{t\in [0,T]}\left\Vert \tilde u_\varepsilon(t)\right\Vert _{\Omega}\leq C, 
\end{equation}
\begin{equation}\label{acotacion3_extension}
 |\tilde u'_\varepsilon|_{p,\Omega,T}\leq C, 
\end{equation}
where the constant $C$ does not depend on $\varepsilon$.
\end{corollary}
\begin{proof}
Using Lemma \ref{extensionCiora} together with the estimate (\ref{acotacion2}) (respectively the estimate (\ref{acotacion5})), we obtain the first estimate in (\ref{acotacion1_extension}) (respectively the estimate (\ref{acotacion2_extension})).

By the Sobolev injection Theorem, if $N=2$ we have the continuous embedding $H_0^1(\Omega)\subset L^p(\Omega)$ and if $N>2$ we have the continuous embedding $H_0^{1}(\Omega)\subset L^{2N/(N-2)}(\Omega)$ which, jointly with the assumption (\ref{assumption_p}), yields the continuous embedding $H_0^1(\Omega)\subset L^p(\Omega)$. 

Therefore, the continuous embedding $H_0^1(\Omega)\subset L^p(\Omega)$ implies that using Lemma \ref{extensionCiora} together with the estimate (\ref{acotacion2}), we can deduce the second estimate in (\ref{acotacion1_extension}).

Finally, using the continuous embedding $H_0^1(\Omega)\subset L^p(\Omega)$, Lemma \ref{extensionCiora} and the estimate (\ref{acotacion8}), we can deduce the estimate (\ref{acotacion3_extension}).
\end{proof}

\section{A compactness result}\label{S5}
In this section, we obtain some compactness results about the behavior of the sequence $\tilde u_\varepsilon$ satisfying the {\it a priori} estimates given in Corollary \ref{estimates_extension}.

Due to the periodicity of the domain $\Omega_\varepsilon$, from Theorem 2.6 in Cioranescu and Donato \cite{CioraDonato} one has, for $\varepsilon \to 0$, that 
\begin{equation}\label{convergence_chi}
\chi_{\Omega_\varepsilon}\stackrel{\tt
*}\rightharpoonup {|Y^*|\over |Y|} \quad \textrm{weakly-star in}\
L^\infty(\Omega),
\end{equation}
where the limit is the proportion of the material in the cell $Y$.

Let $\xi_\varepsilon$ be the gradient of $u_\varepsilon$ in $\Omega_\varepsilon\times (0,T)$ and let us denote by $\tilde \xi_\varepsilon$ its extension with zero to the whole of $\Omega \times (0,T)$, i.e.
\begin{equation}\label{definition_tildexi}
\tilde \xi_\varepsilon=\left\{
\begin{array}{l}
\xi_\varepsilon \quad \text{in }\Omega_\varepsilon \times (0,T),\\
 0 \quad \text{in }(\Omega\setminus \overline{\Omega_\varepsilon})\times (0,T).
 \end{array}\right.
\end{equation}

\begin{proposition}\label{Propo_convergence}
Under the assumptions in Lemma \ref{estimates3}, there exists a function $u\in L^2(0,T;H_0^1(\Omega))\cap L^p(0,T;L^p(\Omega))$ ($u$ will be the unique solution of the limit system (\ref{limit_problem})) and a function $\xi\in L^2(0,T;L^2(\Omega))$ such that for all $T>0,$
\begin{eqnarray}
\label{continuity1} &\tilde{u}_\varepsilon(t)\rightharpoonup
u(t) &\textrm{weakly in}\ H_0^1(\Omega),\quad \forall t\in[0,T],
\\
\label{converge_initial_data}
&\tilde{u}_{\varepsilon}(t)\rightarrow
u(t)&\quad \text{strongly in }L^2(\Omega),\quad \forall t\in[0,T],
\\
\label{converge_f_ae}
&f(\tilde u_{\varepsilon}(t))\rightarrow f(u(t))& \quad \text{strongly in} \quad L^{p'}(\Omega),\quad \forall t\in[0,T],
\\
\label{converge_g_ae}
&g(\tilde u_{\varepsilon}(t))\rightarrow g(u(t))& \quad \text{strongly in} \quad L^{p'}(\Omega),\quad \forall t\in[0,T],
\\
\label{converge_gradiente}
&\tilde \xi_\varepsilon\rightharpoonup \xi& \quad \text{weakly in} \quad L^2(0,T,L^2(\Omega)),
\end{eqnarray}
where $\tilde \xi_\varepsilon$ is given by (\ref{definition_tildexi}).

Moreover, if we suppose that there exists a constant $l>0$ such that
\begin{equation}\label{additional_assumption_g}
\left( g(s)-g(r)\right) \left( s-r\right) \leq l\left(
s-r\right) ^{2},\quad \forall s,r\in\mathbb{R},
\end{equation}
then
\begin{equation}\label{converge_g2_ae}
g(\tilde u_{\varepsilon}(t))\rightharpoonup g(u(t)) \quad \text{weakly in} \quad W_0^{1,p'}(\Omega),\quad \forall t\in[0,T].
\end{equation}
\end{proposition}
\begin{proof}
By (\ref{acotacion1_extension}), we see that the sequence
$\{\tilde u_\varepsilon\}$ is bounded in $L^2(0,T;H_0^1(\Omega))\cap
 L^p(0,T;L^p(\Omega))$, for all $T>0.$ Let us fix $T>0$. Then, there exists a subsequence $\{\tilde u_{\varepsilon'}\}\subset \{\tilde u_\varepsilon\}$ and function $u\in L^2(0,T;H_0^1(\Omega))\cap L^p(0,T;L^p(\Omega))$ such that
\begin{eqnarray}
\label{continuity1_Rellich} &\tilde{u}_{\varepsilon'}\rightharpoonup
u &\textrm{weakly in}\ L^2(0,T;H_0^1(\Omega)),
\\
\label{continuity2} &\tilde{u}_{\varepsilon'}\rightharpoonup
u &\textrm{weakly in}\ L^p(0,T;L^p(\Omega)).
\end{eqnarray}
By the estimate (\ref{acotacion2_extension}), for each $t\in [0,T]$, we have that $\{\tilde u_\varepsilon(t)\}$ is bounded in $H_0^1(\Omega)$, and since we have (\ref{continuity1_Rellich}), we can deduce
\begin{eqnarray*}
 &\tilde{u}_{\varepsilon'}(t)\rightharpoonup
u (t)&\textrm{weakly in}\ H_0^1(\Omega),\quad \forall t\in[0,T].
\end{eqnarray*}
Now, we analyze the convergence for the nonlinear term $f$. 
By the estimate (\ref{acotacion3_extension}), we see that the sequence
$\{\tilde u'_\varepsilon\}$ is bounded in $L^p(0,T;L^p(\Omega))$, for all $T>0.$
Then, we have that $\tilde u_\varepsilon(t):[0,T]\longrightarrow L^p(\Omega)$ is an equicontinuous family of functions.

By Rellich-Kondrachov Theorem and (\ref{assumption_p}), for $p\ge2$, if $N=2$ we have the compact embedding $H_0^1(\Omega)\subset L^p(\Omega)$ and if $N>2$, using that $p\leq 2N/(N-2)$, we also have the compact embedding $H_0^1(\Omega)\subset L^p(\Omega)$.

Since, for each $t\in[0,T]$, we have that $\{\tilde u_\varepsilon(t)\}$ is bounded in $H_0^1(\Omega)$, the compact embedding $H_0^1(\Omega)\subset L^p(\Omega)$, implies that it is precompact in $L^p(\Omega)$.

Then, applying the Ascoli-Arzelà Theorem, we deduce that $\{\tilde u_\varepsilon(t)\}$ is a precompact sequence in $C([0,T];L^p(\Omega))$. Hence, since we have (\ref{continuity2}), we can deduce that
\begin{eqnarray}\label{converge_strongly_p}
 &\tilde{u}_{\varepsilon'}\rightarrow
u &\textrm{strongly in}\ C([0,T];L^p(\Omega)).
\end{eqnarray}
Thanks to (\ref{hipo_consecuencia}), applying Theorem 2.4 in Conca {\it et al.} \cite{Conca} for $G(x,v)=f(v)$, $t=p'$ and $r=p$, we have that the map $v\in L^p(\Omega)\mapsto f(v)\in L^{p'}(\Omega)$ is continuous in the strong topologies. Then, taking into account (\ref{converge_strongly_p}), we get 
$$f(\tilde u_{\varepsilon'}(t))\rightarrow f(u(t)) \quad \text{strongly in} \quad L^{p'}(\Omega)\quad \forall t\in[0,T].$$

Similarly, we analyze the convergence for the nonlinear term $g$ and, we can deduce
\begin{equation}\label{g_fuerte}
g(\tilde u_{\varepsilon'}(t))\rightarrow g(u(t)) \quad \text{strongly in} \quad L^{p'}(\Omega),\quad \forall t\in[0,T].
\end{equation}

In particular, from (\ref{converge_strongly_p}), we have
$$\tilde{u}_{\varepsilon'}(t)\rightarrow
u(t)\quad \text{strongly in }L^p(\Omega), \quad \forall t\in[0,T],\quad \forall p\ge2.$$

To prove (\ref{converge_g2_ae}), let us first note that there exists $C>$ such that
\begin{equation}\label{acotacion_Lqprime}
|\nabla g(\tilde u_\varepsilon(t))|_{p',\Omega}\leq C.
\end{equation}
We observe that under the condition (\ref{additional_assumption_g}), we have that
$$g'(s)\leq l,\quad \forall s\in \mathbb{R}.$$
Then, from the estimate (\ref{acotacion2_extension}), we get
$$\int_{\Omega}\left|{\partial g \over \partial x_i}(\tilde u_\varepsilon(t))\right|^{p'}dx\leq l^{p'}\int_{\Omega}\left|{\partial \tilde u_\varepsilon(t) \over \partial x_i}\right|^{p'}dx\leq C\int_{\Omega}|\nabla \tilde u_\varepsilon(t)|^{p'}dx\leq C\int_{\Omega}|\nabla \tilde u_\varepsilon(t)|^{2}dx\leq C,$$
and we have proved (\ref{acotacion_Lqprime}). Then, from (\ref{g_fuerte}) and (\ref{acotacion_Lqprime}), we can deduce
$$g(\tilde u_{\varepsilon'}(t))\rightharpoonup g(u(t)) \quad \text{weakly in} \quad W_0^{1,p'}(\Omega),\quad \forall t\in[0,T].$$

Finally, from the estimate (\ref{acotacion2}) and (\ref{definition_tildexi}), we have $|\tilde \xi_\varepsilon|_{\Omega,T}\leq C$, and hence, up a sequence, there exists $\xi\in L^2(0,T,L^2(\Omega))$ such that $\tilde \xi_{\varepsilon''} \rightharpoonup \xi$ weakly in $L^2(0,T;L^2(\Omega))$.
  
 By the uniqueness of solution of the limit problem (\ref{limit_problem}), we deduce that the above convergences hold for the whole sequence and therefore, by the arbitrariness of $T>0$, all the convergences are satisfied, as we wanted to prove.

\end{proof}

\section{Homogenized model}\label{S6}
In this section, we identify the homogenized model.
\begin{theorem}\label{Main}
Assume the assumptions in Proposition \ref{Propo_convergence}. Let $(u_\varepsilon, \psi_\varepsilon)$ be the unique solution of the problem (\ref{PDE}). Then, as $\varepsilon\to 0$, we have
$$\tilde u_\varepsilon(t) \to u(t) \quad \text{strongly in } L^2(\Omega),\quad \forall t\in[0,T],$$ 
where $\tilde \cdot$ denotes the extension to $\Omega\times (0,T)$ and $u$ is the unique solution of the following problem
 \begin{equation}\label{limit_problem}
\left\{
\begin{array}{l}
\displaystyle \left({|Y^*|\over |Y|}+{|\partial F| \over |Y|} \right)\displaystyle\frac{\partial u}{\partial t}-\sum_{i,j=1}^Nq_{i,j}{\partial^2u \over \partial x_i \partial x_j}+ {|Y^*|\over |Y|}\left(\kappa u+f(u)\right)+{|\partial F| \over |Y|}g(u)
\\[2ex]
=\displaystyle {|Y^*|\over |Y|}h(x,t)+{|\partial F| \over |Y|}\rho(x,t),   \text{\ in }\;\Omega\times(0,T) ,\\[2ex]
u(x,0)  =  u_{0}(x),  \text{\ for }\;x\in\Omega,\\[2ex]
u= 0,  \text{\ on }
\;\partial \Omega\times( 0,T).
\end{array}
\right.
\end{equation}
The homogenized matrix $Q=((q_{i,j}))$, which is constant and positive-definite, is given by 
\begin{equation}\label{matrix}
q_{i,j}={|Y^*|\over |Y|}\delta_{i,j}-{1\over |Y|}\int_{Y^*}{\partial \eta_j \over \partial y_i}dy,
\end{equation}
where the functions $\eta_j$ are solutions of the system
\begin{equation}\label{system_eta}
\left\{
\begin{array}{l}
\displaystyle -\Delta \eta_j=0,   \text{\ in }Y^*,\\[2ex]
\partial(\eta_j-y_j)/\partial n =0,  \text{\ on }\partial F,\\[2ex]
\eta_j  \text{\ is }Y-\text{periodic},
\end{array}
\right.
\end{equation}
where $y_j$ are local coordinates in $Y^*$.
\end{theorem}

\begin{proof}
We multiply system (\ref{PDE}) by a test function $v\in \mathcal{D}(\Omega)$, and integrating by parts, we have
\begin{eqnarray*}
\dfrac{d}{dt}\left(\int_{\Omega}\chi_{\Omega_\varepsilon}\tilde u_\varepsilon(t)vdx\right)+\varepsilon\,\dfrac{d}{dt}\left(\int_{\partial F_\varepsilon}
\gamma_{0}(u_\varepsilon(t))vd\sigma(x)\right)+\int_{\Omega}\tilde \xi_\varepsilon \nabla vdx+\kappa \int_{\Omega}\chi_{\Omega_\varepsilon}\tilde u_\varepsilon(t)v dx\\[2ex]
+\int_{\Omega}\chi_{\Omega_\varepsilon} f(\tilde u_\varepsilon(t))vdx
 +\varepsilon\, \int_{\partial F_\varepsilon}
g(\gamma_{0}(u_\varepsilon(t)))vd\sigma(x)
 =\int_{\Omega} \chi_{\Omega_\varepsilon} h(t)vdx
 +\varepsilon\int_{\partial F_\varepsilon}\rho(t)vd\sigma(x), 
 \end{eqnarray*}
 in $\mathcal{D}'(0,T)$.
 
 We consider $\varphi\in C_c^1([0,T])$ such that $\varphi(T)=0$ and $\varphi(0)\ne 0$. Multiplying by $\varphi$ and integrating between $0$ and $T$, we have
 \begin{eqnarray}\label{system1}\nonumber
-\varphi(0)\left(\int_{\Omega}\chi_{\Omega_\varepsilon}\tilde u_\varepsilon(0)vdx\right)-\int_0^T\dfrac{d}{dt}\varphi(t)\left(\int_{\Omega}\chi_{\Omega_\varepsilon}\tilde u_\varepsilon(t)vdx\right)dt\\[2ex]\nonumber
-\varepsilon\varphi(0)\left(\int_{\partial F_\varepsilon}
\gamma_{0}(u_\varepsilon(0))vd\sigma(x)\right)
-\varepsilon\int_0^T \dfrac{d}{dt}\varphi(t)\left(\int_{\partial F_\varepsilon}
\gamma_{0}(u_\varepsilon(t))vd\sigma(x)\right)dt\\[2ex]
+\int_0^T\varphi(t)\int_{\Omega}\tilde \xi_\varepsilon\nabla vdxdt+\kappa \int_0^T\varphi(t)\int_{\Omega}\chi_{\Omega_\varepsilon}\tilde u_\varepsilon(t)v dxdt\\[2ex]\nonumber
+\int_0^T \varphi(t)\int_{\Omega}\chi_{\Omega_\varepsilon}f(\tilde u_\varepsilon(t))vdxdt
 +\varepsilon \int_0^T \varphi(t)\int_{\partial F_\varepsilon}
g(\gamma_{0}(u_\varepsilon(t)))vd\sigma(x)dt\\[2ex]
 =\int_0^T \varphi(t)\int_{\Omega}\chi_{\Omega_\varepsilon} h(t)vdxdt
 +\varepsilon \int_0^T \varphi(t)\int_{\partial F_\varepsilon}\rho(t)vd\sigma(x)dt.\nonumber
 \end{eqnarray}

For the sake of clarity, we split the proof in three parts. Firstly, for the integrals on $\Omega$ we only require to use Proposition \ref{Propo_convergence} and the convergence (\ref{convergence_chi}), secondly for the integrals on the boundary of the holes we make use of a convergence result based on a technique introduced by Vanninathan \cite{Vanni}. Finally, we pass to the limit, as $\varepsilon \to 0$, in (\ref{system1}).

{\bf Step 1}. Passing to the limit, as $\varepsilon\to 0$, in the integrals on $\Omega$:

From (\ref{converge_initial_data})-(\ref{converge_f_ae}) and (\ref{convergence_chi}), we have respectively, for $\varepsilon \to 0$, 
\begin{equation*}
\int_{\Omega}\chi_{\Omega_\varepsilon}\tilde u_\varepsilon(t)v dx\to {|Y^*|\over |Y|}\int_{\Omega}u(t)v dx,\quad \forall v \in \mathcal{D}(\Omega),
\end{equation*}
and 
\begin{equation*}
\int_{\Omega}\chi_{\Omega_\varepsilon}f(\tilde u_\varepsilon(t)) vdx\to {|Y^*|\over |Y|}\int_{\Omega}f(u(t))vdx,\quad \forall v\in \mathcal{D}(\Omega),
\end{equation*}
which integrating in time and using Lebesgue's Dominated Convergence Theorem, gives
\begin{equation*}
\int_0^T{d \over dt}\varphi(t)\left(\int_{\Omega}\chi_{\Omega_\varepsilon}\tilde u_\varepsilon(t) vdx\right)dt\to {|Y^*|\over |Y|}\int_0^T{d \over dt}\varphi(t)\left(\int_{\Omega}u(t)vdx\right)dt,
\end{equation*}
\begin{equation*}
\kappa\int_0^T\varphi(t)\int_{\Omega}\chi_{\Omega_\varepsilon}\tilde u_\varepsilon(t) vdxdt\to \kappa{|Y^*|\over |Y|}\int_0^T\varphi(t)\int_{\Omega}u(t)vdxdt,
\end{equation*}
and 
\begin{equation*}
\int_0^T\varphi(t)\int_{\Omega}\chi_{\Omega_\varepsilon}f(\tilde u_\varepsilon(t)) vdxdt\to {|Y^*|\over |Y|}\int_0^T\varphi(t)\int_{\Omega}f(u(t))vdxdt.
\end{equation*}
By (\ref{converge_initial_data}) and (\ref{convergence_chi}), we have
\begin{equation*}
\varphi(0)\left(\int_{\Omega}\chi_{\Omega_\varepsilon}\tilde u_\varepsilon(0)vdx\right)\to \varphi(0){|Y^*|\over |Y|}\int_{\Omega} u(0)vdx, \quad \forall v\in\mathcal{D}(\Omega).
\end{equation*}
By the assumption (\ref{hyp 0'}), (\ref{convergence_chi}) and using Lebesgue's Dominated Convergence Theorem, we get
$$\int_0^T \varphi(t)\int_{\Omega}\chi_{\Omega_\varepsilon} h(t)vdxdt \to {|Y^*|\over |Y|}\int_0^T\varphi(t)\int_{\Omega}h(t)vdxdt.$$

On the other hand, using (\ref{converge_gradiente}), we obtain, for $\varepsilon\to 0$
$$\int_0^T\varphi(t)\int_{\Omega}\tilde \xi_\varepsilon\nabla vdxdt \to \int_0^T\varphi(t)\int_{\Omega}\xi\nabla vdxdt.$$
{\bf Step 2}. Passing to the limit, as $\varepsilon\to 0$, in the surface integrals on the boundary of the holes:

We make use of the technique introduced by Vanninathan \cite{Vanni} for the Steklov problem which transforms surface integrals into volume integrals. This technique was already used as a main tool to homogenize the non homogeneous Neumann problem for the elliptic case by Cioranescu and Donato \cite{Ciora2}. 

By Definition 3.2 in Cioranescu and Donato \cite{Ciora2}, let us introduce, for any $h\in L^{s'}(\partial F)$, $1\leq s'\leq \infty$, the linear form $\mu_{h}^\varepsilon$ on $W_0^{1,s}(\Omega)$ defined by
$$\langle \mu_{h}^\varepsilon,\varphi \rangle=\varepsilon \int_{\partial F_\varepsilon} h\left(x\over \varepsilon \right)\varphi d\sigma,\quad \forall \varphi\in W_0^{1,s}(\Omega),$$
with $1/s+1/s'=1$. It is proved in Lemma 3.3 in Cioranescu and Donato \cite{Ciora2} that
\begin{equation}\label{convergence_mu}
\mu_{h}^\varepsilon \to \mu_{h}\quad \text{strongly in }(W_0^{1,s}(\Omega))',
\end{equation}
where $\langle \mu_{h},\varphi \rangle=\mu_{h}\displaystyle\int_{\Omega}\varphi dx$, with
$$\mu_{h}={1\over |Y|}\int_{\partial F}h(y)d\sigma.$$
In the particular case in which $h\in L^{\infty}(\partial F)$ or even when $h$ is constant, we have
\begin{equation*}
\mu_{h}^\varepsilon \to \mu_{h}\quad \text{strongly in }W^{-1,\infty}(\Omega).
\end{equation*}
In what follows, we shall denote by $\mu_1^\varepsilon$ the above introduced measure in the particular case in which $h=1$. Notice that in this case $\mu_{h}$ becomes $\mu_1=|\partial F|/|Y|$.

Observe that using Corollary 4.2 in Cioranescu {\it et al.} \cite{Ciora3} with (\ref{continuity1}), we can deduce, for $\varepsilon \to 0$, 
\begin{equation*}
\varepsilon \int_{\partial F_\varepsilon} \gamma_0(u_\varepsilon(t))v d\sigma(x)=\langle \mu_{1}^\varepsilon,\tilde u_{\varepsilon {|_{\Omega_\varepsilon}}}(t)v \rangle\to \mu_1\int_{\Omega}u(t)v dx={|\partial F| \over |Y|}\int_{\Omega}u(t)v dx,
\end{equation*}
for all $v\in \mathcal{D}(\Omega)$, which integrating in time and using Lebesgue's Dominated Convergence Theorem, gives
$$\varepsilon \int_0^T{d \over dt}\varphi(t)\left(\int_{\partial F_\varepsilon} \gamma_0(u_\varepsilon(t))v d\sigma(x)\right) dt \to {|\partial F| \over |Y|}\int_0^T{d \over dt}\varphi(t)\left(\int_{\Omega}u(t)v dx\right)dt.$$
Moreover, using Corollary 4.2 in Cioranescu {\it et al.} \cite{Ciora3} with (\ref{continuity1}), we can deduce, for $\varepsilon \to 0$, 
\begin{equation*}
\varepsilon\int_{\partial F_\varepsilon}
\gamma_{0}(u_\varepsilon(0))vd\sigma(x)=\langle \mu_{1}^\varepsilon,\tilde u_{\varepsilon {|_{\Omega_\varepsilon}}}(0)v \rangle\to \mu_1\int_{\Omega}u(0)v dx={|\partial F| \over |Y|}\int_{\Omega}u(0)v dx, \quad \forall v\in \mathcal{D}(\Omega).
\end{equation*}
On the other hand, note that using (\ref{convergence_mu}) with $s=2$, taking into account (\ref{hyp 0'}) and by Lebesgue's Dominated Convergence Theorem, we can deduce, for $\varepsilon \to 0$, 
$$\varepsilon \int_0^T \varphi(t)\int_{\partial F_\varepsilon}\rho(t)vd\sigma(x)dt=\int_0^T\varphi(t)\langle \mu_{1}^\varepsilon,\rho(t)v \rangle dt\to {|\partial F| \over |Y|}\int_0^T \varphi(t)\int_{\Omega}\rho(t)vdxdt.$$

From (\ref{converge_g2_ae}) and (\ref{convergence_mu}), with $s=p'$, we conclude
$$\varepsilon \int_{\partial F_\varepsilon}
g(\gamma_{0}(u_\varepsilon(t)))vd\sigma(x)=\langle \mu_{1}^\varepsilon,g(\tilde u_{\varepsilon}(t))v \rangle\to {|\partial F| \over |Y|}\int_{\Omega}g(u(t))vdx,$$
for all $v\in \mathcal{D}(\Omega)$, which integrating in time and using Lebesgue's Dominated Convergence Theorem, gives
$$\varepsilon \int_0^T\varphi(t)\int_{\partial F_\varepsilon}
g(\gamma_{0}(u_\varepsilon(t)))vd\sigma(x)dt \to {|\partial F| \over |Y|}\int_0^T\varphi(t)\int_{\Omega}g(u(t))vdxdt.$$

{\bf Step 3}. Passing to the limit, as $\varepsilon\to 0$, in (\ref{system1}):

All the terms in (\ref{system1}) pass to the limit, as $\varepsilon \to 0$, and therefore taking into account the previous steps, we get
 \begin{eqnarray*}
-\varphi(0)\left({|Y^*|\over |Y|}+{|\partial F| \over |Y|} \right)\left(\int_{\Omega} u(0)vdx\right)-\left({|Y^*|\over |Y|}+{|\partial F| \over |Y|} \right)\int_0^T\dfrac{d}{dt}\varphi(t)\left(\int_{\Omega} u(t)vdx\right)dt\\[2ex]\nonumber
+\int_0^T\varphi(t)\int_{\Omega}\xi\nabla vdxdt+\kappa {|Y^*|\over |Y|}\int_0^T\varphi(t)\int_{\Omega} u(t)v dxdt\\[2ex]\nonumber
+{|Y^*|\over |Y|}\int_0^T \varphi(t)\int_{\Omega}f( u(t))vdxdt
 +{|\partial F| \over |Y|} \int_0^T \varphi(t)\int_{\Omega}
g(u(t))vdxdt\\[2ex]
 ={|Y^*|\over |Y|}\int_0^T \varphi(t)\int_{\Omega} h(t)vdxdt
 +{|\partial F| \over |Y|} \int_0^T \varphi(t)\int_{\Omega}\rho(t)vdxdt.\nonumber
 \end{eqnarray*}
 Hence, $\xi$ verifies
 \begin{equation}\label{equation_xi}
 \left({|Y^*|\over |Y|}+{|\partial F| \over |Y|} \right)\displaystyle\frac{\partial u}{\partial t}-{\rm div}\xi+ {|Y^*|\over |Y|}\left(\kappa u+f(u)\right)+{|\partial F| \over |Y|}g(u)= {|Y^*|\over |Y|}h+{|\partial F| \over |Y|}\rho, \quad \text{in }\Omega\times (0,T).
 \end{equation}
  
It remains now to identify $\xi$. The proof is standard, so we omit it. Following, for example, the proof of Theorem 4.7 in Cioranescu and Donato \cite{Ciora2}, we conclude that 
 \begin{equation}\label{elliptical_case}
 \xi=Q\nabla u,\quad \text{in } \Omega\times (0,T),
 \end{equation}
 where $Q=((q_{i,j}))$ is given by (\ref{matrix}). Then, taking into account (\ref{elliptical_case}) in (\ref{equation_xi}), we have the homogenized model (\ref{limit_problem}).
 
\end{proof}
 
 \begin{definition}
A weak solution of (\ref{limit_problem}) is any function $u$, satisfying
\begin{equation*}
u\in {C}([0,T];L^{2}\left(  \Omega\right)  ),\quad \text{for all }T>0,
\end{equation*}
\begin{equation*}
u\in L^{2}(0,T;H_{0}^{1}\left(  \Omega\right)  )\cap L^{p}(0
,T;L^{p}\left(  \Omega\right)  ),\quad \text{for all }T>0,
\end{equation*}
\begin{equation*}
\!\!\!\!\!\!\left\{
\begin{array}{l}
\displaystyle \left({|Y^*|\over |Y|}+{|\partial F| \over |Y|} \right) \dfrac{d}{dt}(u(t),v)+(Q\nabla u(t),\nabla v)+{|Y^*|\over |Y|}\kappa(u(t),v)+{|Y^*|\over |Y|}(f(u(t)),v)+{|\partial F| \over |Y|}(g(u(t)),v)\\[2ex]
\displaystyle={|Y^*|\over |Y|}(h(t),v)+{|\partial F| \over |Y|}(\rho(t),v),\\[2ex]
 \hbox{in $\mathcal{D}'(0,T)$, for all $v\in H_0^1(\Omega)\cap L^p(\Omega)$,}
\end{array}
\right.
\end{equation*}
\begin{equation*}
 u(0)=u_0.
\end{equation*}
\end{definition}

\begin{remark}
Applying a slight modification of Theorem 1.4, Chapter 2 in Lions \cite{Lions}, we obtain that the problem (\ref{limit_problem}) has a unique solution.
\end{remark}

\subsection*{Acknowledgments}
The author has been supported by Junta de Andaluc\'ia (Spain), Proyecto de Excelencia P12-FQM-2466.

\end{document}